\documentclass[12pt]{article}
\pdfoutput=1
\usepackage{graphicx,pict2e,latexsym}
\usepackage[square,numbers]{natbib}
\newtheorem{theorem}{Theorem}
\newtheorem{lemma}{Lemma}
\newtheorem{corollary}{Corollary}

\begin{document}
\bibliographystyle{plainnat}
\pagestyle{plain}

\title{\Large \bf On the Power of Likelihood Ratio Tests in Dimension-Restricted Submodels}

\author{Michael W. Trosset\thanks{Department of Statistics, Indiana University.  
E-mail: {\tt mtrosset@indiana.edu}} \and
Mingyue Gao\thanks{Department of Applied Mathematics \& Statistics, Johns Hopkins University. 
E-mail: {\tt mygao90@gmail.com}} \and
Carey E. Priebe\thanks{Department of Applied Mathematics \& Statistics, Johns Hopkins University. 
E-mail: {\tt cep@jhu.edu}} \and
}

\date{\today}

\maketitle

\newpage

\begin{abstract}
Likelihood ratio tests are widely used to test statistical hypotheses about parametric families of probability distributions.
If interest is restricted to a subfamily of distributions, then it is natural to inquire if the restricted LRT is superior to the unrestricted LRT.  
Marden's general LRT conjecture posits that any restriction placed on the alternative hypothesis will increase power.
The only published counterexample to this conjecture is rather technical and involves a restriction that maintains the dimension of the alternative.
We formulate the dimension-restricted LRT conjecture, which posits that any restriction that replaces a parametric family with a subfamily of lower dimension will increase power.
Under standard regularity conditions, we then demonstrate that the restricted LRT is asymptotically more powerful than the unrestricted LRT for local alternatives.
Remarkably, however, even the dimension-restricted LRT conjecture fails in the case of finite samples.
Our counterexamples involve subfamilies of multinomial distributions.  In particular, our study of
the Hardy-Weinberg subfamily of trinomial distributions provides a simple and elegant demonstration that restrictions may not increase power.
\end{abstract}


\tableofcontents

\newpage


\section{Introduction}
\label{intro}

We compare restricted and unrestricted likelihood ratio tests
in situations where the restriction decreases the dimension of the alternative.  The issues that concern us are motivated by an elementary example.

\subparagraph{Basic Example}  Suppose that $X=(X_1,X_2)$ has a bivariate normal distribution with mean vector $\theta = (\theta_1,\theta_2)$ and identity covariance matrix, in which case the parametric family of possible probability distributions (the model) is $2$-dimensional.  Let $\vec{0}$ denote the origin in $\Re^2$ and consider testing the simple null hypothesis 
$H_0: \theta = \vec{0}$ against the $2$-dimensional composite alternative hypothesis $H_A: \theta \neq \vec{0}$.  Under this model, 
the likelihood function is
\[
L_x(\theta) = \frac{1}{2 \pi} \exp \left( -\frac{1}{2} \left\| x-\theta \right\|^2 \right),
\]
the (unrestricted) maximum likelihood estimate (MLE) of $\theta$ is $\hat{\theta}=x$, 
and the (unrestricted) likelihood ratio test (LRT) rejects $H_0$
if and only if
\[
\Lambda_2(x) = -2 \log \frac{L_x(\vec{0})}{L_x(\hat{\theta})} = -2 \log
\frac{\exp \left( -\frac{1}{2} \| x-0 \|^2 \right)}{\exp \left( -\frac{1}{2} \| x-x \|^2 \right)} =  \| x \|^2 = x_1^2 + x_2^2
\]
is sufficiently large.  Because the random variable $X_1^2+X_2^2 \sim \chi^2_2$, a (central) chi-squared distribution with $2$ degrees of freedom, the unrestricted LRT rejects $H_0$ at significance level $\alpha$ if and only if
\[
\phi_2(x) = \| x \|^2 > c_{2,\alpha},
\]
where $c_{2,\alpha}$ is the $1-\alpha$ quantile of $\chi^2_2$.

Suppose that we know that $\theta = (\theta_1,\theta_2) = (\theta_1,0)$.  Restricting attention to this $1$-dimensional submodel, we can write the null hypothesis as $H_0: \theta_1 = 0$ and the alternative hypothesis as $H_A: \theta_1 \neq 0$.
Under this submodel,
the likelihood function is
\[
L_x(\theta) = \frac{1}{2 \pi} \exp \left( -\frac{1}{2} \left[ \left( x_1-\theta_1 \right)^2 + \left( x_2-0\right)^2 \right] \right),
\] 
the (restricted) maximum likelihood estimate (MLE) of $\theta$ is $\tilde{\theta}=(x_1,0)$, 
and the (restricted) likelihood ratio test (LRT) rejects $H_0$
if and only if
\[
\Lambda_1(x) = -2 \log \frac{L_x(\vec{0})}{L_x(\tilde{\theta})} = -2 \log
\frac{\exp \left( -\frac{1}{2} \left[ x_1^2+x_2^2 \right] \right)}{\exp \left( -\frac{1}{2} x_2^2 \right)} =  x_1^2 
\]
is sufficiently large.  Because the random variable $X_1^2 \sim \chi^2_1$,  the restricted LRT rejects $H_0$ at significance level $\alpha$ if and only if
\[
\phi_1(x) = x_1^2 > c_{1,\alpha},
\]
where $c_{1,\alpha}$ is the $1-\alpha$ quantile of $\chi^2_1$.

We now compare the power of $\phi_1$ and $\phi_2$ at alternatives of the form $\theta=(\delta,0)$.  If $X_1 \sim \mbox{Normal}(\delta,1)$ and $X_2 \sim \mbox{Normal}(0,1)$, then 
\begin{eqnarray*}
X_1^2 \sim \chi^2_1(\delta^2) & \mbox{ and } & 
X_1^2+X_2^2 \sim \chi^2_2(\delta^2),
\end{eqnarray*}
where $\chi^2_m(\lambda)$ is the noncentral chi-squared distribution with $m$ degrees of freedom and noncentrality parameter $\lambda$.
The power of $\phi_1$ is
\[
\pi_\alpha \left( \delta^2; 1 \right) =
P \left( \chi^2_1 \left( \delta^2 \right) > c_{1,\alpha} \right),
\]
the power of $\phi_2$ is
\[
\pi_\alpha \left( \delta^2; 2 \right) =
P \left( \chi^2_2 \left( \delta^2 \right) > c_{2,\alpha} \right),
\]
and the following lemma implies that
\[
\pi_\alpha \left( \delta^2; 1 \right) >
\pi_\alpha \left( \delta^2; 2 \right).
\]
\hfill $\Box$

\bigskip
\noindent
{\bf Lemma} {$\mathbf \chi^2$} (Das Gupta and Perlman
\citep{dasgupta&perlman:1974}):
{\em Fix $\alpha \in (0,1)$ and $\lambda>0$.
For $m \geq 1$, define $c_{m,\alpha}$ by
$P(\chi^2_m > c_{m,\alpha}) = \alpha$.
Then
\[
\pi_\alpha ( \lambda; m ) =
P \left( \chi^2_m ( \lambda ) > c_{m,\alpha} \right)
\]
is strictly decreasing in $m$.}

\bigskip

The essential property of the Basic Example is that the restricted LRT is more powerful on the submodel than the unrestricted LRT.
Numerous studies have explored the extent to which this property can be generalized.  The most studied problem is that of order-restricted inference on normal means, exemplified in the Basic Example by restricting attention to  the $2$-dimensional submodel defined by $\theta_1 \leq \theta_2$.
The most general result of this type appears to be that of
Praestgaard \citep{praestgaard:2012}:
\begin{quote}
``Let $X_1,X_2,\ldots,X_n$ be normally distributed random variables with means
$\vec{\mu}=(\mu_1,\mu_2,\ldots,\mu_n)$ and known positive definite covariance matrix $\Sigma$.  Let ${\mathcal C}$ be a closed convex cone in $\Re^n$ which contains a linear space ${\mathcal L} \subset {\mathcal C}$.  The present note considers the likelihood ratio test for null and alternative hypotheses
\begin{eqnarray*}
H_0: \vec{\mu} \in {\mathcal L} & \mbox{ versus } &
H_1: \vec{\mu} \in {\mathcal C}/{\mathcal L}.
\end{eqnarray*}
We prove that the restricted likelihood ratio test is uniformly more powerful over 
$\vec{\mu} \in {\mathcal C}$ than the omnibus test with
alternative hypothesis $\vec{\mu} \in \Re^n$.''
\end{quote}

Studies of restricted LRTs for order-restricted inference date to the pioneering work of Bartholomew 
\citep{bartholomew:1959a,bartholomew:1959b,bartholomew:1961}
in the late 1950s and early 1960s.
By the 1980s, it was imagined that the power superiority of restricted LRTs might be a universal phenomenon.
Al-Rawwash's 1986 Ph.D. dissertation 
\citep{alrawwash:1986}
stated a general LRT conjecture (``the more restrictions which are put on the alternative space, the higher the power of the L.R.T.''), attributing it to a 1982 NSF proposal submitted by his advisor, J. Marden; however,
Al-Rawwash \citep[Chapter 6]{alrawwash:1986}
only studied the Basic Example with conic submodels.

As late as 1992, Tsai \citep{tsai:1992} was able to state that
``The long time conjecture of the power superiority of the restricted LRT to
its unrestricted version in the entire parameter space of alternatives for the general setting is of considerably analytic difficulty and lack of the
definitive results.''
In 2003, however, Abu-Dayyeh, Al-Jararha, and Madan
\citep{abudayyeh&etal:2003} constructed a counterexample using the Basic Example with nonconic submodels of the form $\theta \in [-k,k]^2$.
Surprisingly, smaller values of $k$ give more restricted alternatives,
but {\em not}\/ uniformly greater power.

In our presentation of the Basic Example, the power superiority of the restricted LRT appears to be a consequence of the fact that the chi-squared distribution of the restricted test statistic has fewer degrees of freedom than the chi-squared distribution of the unrestricted test statistic.  This fact, in turn, is a consequence of the fact that the dimension of the submodel is smaller than the dimension of the model.  Noting that the dimension of the submodels $\theta \in [-k,k]^2$ is the same as the dimension of the model,
it is inviting to modify the general LRT conjecture by speculating that it holds if the submodel has lower dimension than the model.
We will refer to this special case of Marden's general LRT conjecture as the dimension-restricted LRT conjecture.
One would expect the dimension-restricted LRT conjecture to hold at least asymptotically because (under suitable regularity conditions) LRT statistics have asymptotic chi-squared distributions.  In Section \ref{asymptotic} we exploit that fact and demonstrate that the dimension-restricted LRT conjecture {\em does}\/ hold asymptotically.  The more interesting question is whether or not it holds for finite samples.

Despite considerable interest in Marden's general LRT conjecture, we are not aware of any previous statements of the more plausible dimension-restricted LRT conjecture.  We demonstrate in Section \ref{HW} that the dimension-restricted LRT conjecture does not hold for finite samples.  To construct a counterexample, we abandon normal models and study the $1$-dimensional Hardy-Weinberg submodel of the $2$-dimensional trinomial model.  The implications of this counterexample are considered in Section~\ref{disc}.

\section{Asymptotic Theory}
\label{asymptotic}

We begin by recalling some basic properties of differentiable manifolds,
which will serve as index sets for our statistical models and submodels.
Let $M$ denote a completely separable Hausdorff space.  
(For our purposes, it will suffice to assume that $M$ is a subset of Euclidean space.)
Let $U \subseteq M$ and $V \subseteq \Re^k$ denote open sets.
If $\varphi : U \rightarrow V$ is a homeomorphism,
then $\varphi(u) = (x_1(u),\ldots,x_k(u))$ defines a coordinate system on $U$.
The $x_i$ are the coordinate functions and $\varphi^{-1}$ is a parametrization of $U$.  The pair $(U,\varphi)$ is a chart.
An atlas on $M$ is a collection of charts $\{ (U_a,\varphi_a) \}$ such that the $U_a$ cover $M$.

The set $M$ is a $k$-dimensional topological manifold if and only if it admits an atlas for which each $\varphi_a(U_a)$ is open in $\Re^k$.
It is smooth if and only if the transition maps $\varphi_b \varphi_a^{-1}$ are diffeomorphisms.
A subset $S \subset M$ is a $d$-dimensional embedded submanifold if and only if, for every $p \in S$, there is a chart $(U,\varphi)$ such that $p \in U$
and
\[
\varphi(U \cap S) = \varphi(U) \cap \left( \Re^d \times 
\{ \vec{0} \in \Re^{k-d} \} \right) =
\left\{ y \in \varphi(U) : y_{d+1}= \cdots = y_k=0 \right\}.
\]
We will
assume that the statistical model $\{ P_\theta : \theta \in \Theta \}$
is indexed by a $k$-dimensional manifold $\Theta$ and that the statistical submodel $\{ P_\theta : \theta \in \Psi \}$ is indexed by a $d$-dimensional embedded submanifold $\Psi \subset \Theta$.

Let $X_1,\ldots,X_n$ be independent and identically distributed as $P_\theta$ and suppose that $\Theta_0 \subset \Psi$ is an $\ell$-dimensional manifold.
For testing $H_0 : \theta \in \Theta_0$,
the unrestricted LRT statistic is
\[
\Lambda_{k,n} \left( X_1,\ldots,X_n \right) = -2 \log
\frac{\sup_{\theta \in \Theta_0} \prod_{i=1}^n p_\theta (X_i)}{ 
\sup_{\theta \in \Theta} \prod_{i=1}^n p_\theta (X_i)}
\]
and the restricted LRT statistic is
\[
\Lambda_{d,n} \left( X_1,\ldots,X_n \right) = -2 \log
\frac{\sup_{\theta \in \Theta_0} \prod_{i=1}^n p_\theta (X_i)}{ 
\sup_{\theta \in \Psi} \prod_{i=1}^n p_\theta (X_i)}.
\]
We will compare the power of these LRTs at local alternatives in
$\Psi/\Theta_0$ as $n \rightarrow \infty$.
To do so, we require the asymptotic distributions of $\Lambda_{k,n}$ and $\Lambda_{d,n}$.

Under classical regularity conditions, the asymptotic null distribution of $\Lambda_{k,n}$ is $\chi^2_{k-\ell}$ \citep{wilks:1938,chernoff:1954}.
The same conclusion was drawn by 
van der Vaart \citep[Chapter 16]{vandervaart:1998}, who based his derivation on the convergence of experiments.  Using this approach, ``the main conditions are that the model is differentiable in $\theta$ [more precisely, that the map $\theta \mapsto \sqrt{p_\theta}$ is differentiable in quadratic mean] and that the null hypothesis $\Theta_0$ and the full parameter set $\Theta$ are (locally) equal to linear spaces [i.e., $\Theta_0$ and $\Theta$ are manifolds].''  By the same reasoning, of course, the asymptotic null distribution of $\Lambda_{d,n}$ is $\chi^2_{d-\ell}$.
We might assume any set of conditions that ensure these asymptotic distributions; for greatest generality, we  simply assume that,
for $\vartheta \in \Theta_0$,
\begin{eqnarray}
\lim_{n \rightarrow \infty} P_\vartheta \left( \Lambda_{k,n} > c \right) =
P \left( \chi^2_{k-\ell} > c \right) & \mbox{ and } &
\lim_{n \rightarrow \infty} P_\vartheta \left( \Lambda_{d,n} > c \right) =
P \left( \chi^2_{d-\ell} > c \right).
\label{eq:null}
\end{eqnarray}

Local asymptotic power functions were also derived (in the special case of local asymptotic normality, although the approach is more general) by
van der Vaart \citep[Section 16.4]{vandervaart:1998}.
Local alternatives are alternatives of the form $\theta = \vartheta + h/\sqrt{n}$ for $\vartheta \in \Theta_0$ and $n$ sufficiently large.
In order to compare the restricted and unrestricted LRTs,
we study power at local alternatives $\vartheta + h/\sqrt{n} \in \Psi$.

Let $I_k(\vartheta)$ and $I_d(\vartheta)$ denote the Fisher information matrices at $\vartheta$ for the model and submodel respectively.
Following van der Vaart \citep[Chapter 16]{vandervaart:1998}, 
define sets
\[
H_{n,0} = \sqrt{n} \left( \Theta_0 - \vartheta \right) =
\left\{ \sqrt{n} \left( \vartheta^\prime - \vartheta \right) :
\vartheta^\prime \in \Theta_0 \right\}
\]
and let $H_0$ denote the set of all limits of convergent sequences
$\{ h_n \in H_{n,0} \}$.
Assume that the limits of all convergent subsequences
$\{ h_{n_i} \in H_{n_i,0} \}$ lie in $H_0$, in which case we write
$H_{n,0} \rightarrow H_0$.  In the case of a simple null hypothesis $\Theta_0 = \{ \theta_0 \}$,
\[
H_{n,0} = \sqrt{n} \left( \theta_0 - \theta_0 \right) = 
\{ \vec{0} \in \Re^k \} = H_0.
\]
Then, under suitable regularity conditions, 
\begin{equation}
\lim_{n \rightarrow \infty} P_{\vartheta + h/\sqrt{n}} 
\left( \Lambda_{k,n} > c \right) =
P \left( \chi^2_{k-\ell} \left( \lambda_k \right) > c \right),
\label{eq:altk1}
\end{equation}
with noncentrality parameter
\begin{equation}
\lambda_k = \delta_k^2 =
\inf_{h^\prime \in H_0} \left( h-h^\prime \right)^T 
I_k (\vartheta) \left( h-h^\prime \right),
\label{eq:altk2}
\end{equation}
and
\begin{equation}
\lim_{n \rightarrow \infty} P_{\vartheta + h/\sqrt{n}} 
\left( \Lambda_{d,n} > c \right) =
P \left( \chi^2_{d-\ell} \left( \lambda_d \right) > c \right),
\label{eq:altd1}
\end{equation}
with noncentrality parameter
\begin{equation}
\lambda_d = \delta_d^2 =
\inf_{h^\prime \in H_0} \left( h-h^\prime \right)^T 
I_d (\vartheta) \left( h-h^\prime \right),
\label{eq:altd2}
\end{equation}
Again, for greatest generality, we simply assume that
(\ref{eq:altk1})--(\ref{eq:altd2}) 
hold.

To apply Lemma~$\chi^2$, we will require
\bigskip
\begin{lemma}
If $\Psi$ is an embedded submanifold of $\Theta$, then
$\lambda_k=\lambda_d$.
\label{lm:noncentral}
\end{lemma}

\subparagraph{Proof} 
Because $\varphi \in \Theta_0 \subset \Psi$ and $\Psi$ is an embedded submanifold of $\Theta$, there exists a chart $(U,\varphi)$ such that
$\vartheta \in U$ and 
\[
\varphi(u) = (y_1,\ldots,y_d,0,\ldots,0) \in \Re^k
\]
if $u \in U \cap \Psi$.
Choose $n$ large enough that $\vartheta + h/\sqrt{n} \in U$, then write
(\ref{eq:altk2}) and (\ref{eq:altd2}) in the coordinate system defined by $\varphi$.  We obtain
\begin{eqnarray*}
\lambda_k & = &
\inf_{(z_1,\ldots,z_d,0,\ldots,0) \in H_0} 
\left[ \begin{array}{ccc|ccc}
y_1-z_1 & \ldots & y_d-z_d & 0 & \ldots & 0
\end{array} \right]
\left[ \begin{array}{c|c} 
I_d (\vartheta) & \cdot \\ \hline
\cdot & \cdot
\end{array} \right]
\left[ \begin{array}{c}
y_1-z_1 \\ \vdots \\ y_d-z_d \\ \hline 0 \\ \vdots \\ 0
\end{array} \right] \\
 & = &
\inf_{(z_1,\ldots,z_d) \in H_0} 
\left[ \begin{array}{ccc}
y_1-z_1 & \ldots & y_d-z_d 
\end{array} \right]
I_d (\vartheta) 
\left[ \begin{array}{c}
y_1-z_1 \\ \vdots \\ y_d-z_d 
\end{array} \right] \\
 & = & \lambda_d.
\end{eqnarray*}
\hfill $\Box$

\bigskip
Now we establish the asymptotic power superiority of the restricted LRT.
Given $\alpha \in (0,1)$,
define quantiles $c_{k,\alpha}$ and $c_{d,\alpha}$ by
\[
P \left( \chi^2_{k-\ell} > c_{k,\alpha} \right) = \alpha =
P \left( \chi^2_{d-\ell} > c_{d,\alpha} \right).
\]
The unrestricted LRT rejects $H_0: \theta \in \Theta_0$ if and only if
\[
\phi_k \left( \vec{x}_n \right) = 
\Lambda_{k,n} \left( x_1,\ldots,x_n \right)  > c_{k,\alpha},
\]
and the restricted LRT rejects $H_0: \theta \in \Theta_0$ if and only if
\[
\phi_d \left( \vec{x}_n \right) = 
\Lambda_{d,n} \left( x_1,\ldots,x_n \right)  > c_{d,\alpha}.
\]
We compare the power of $\phi_k$ and $\phi_d$  
at local alternatives $\vartheta + h/\sqrt{n} \in \Psi$.
\bigskip
\begin{theorem}
Fix $\vartheta \in \Theta_0$ and $h$ such that $\vartheta + h/\sqrt{n} \in \Psi$ for $n \geq N$.  If $\Psi$ is an embedded submanifold of $\Theta$,
then
\[
P_{\vartheta + h/\sqrt{n}} \left( 
\Lambda_{d,n}   > c_{d,\alpha}
\right) >
P_{\vartheta + h/\sqrt{n}} \left( 
\Lambda_{k,n}   > c_{k,\alpha}
\right)
\]
for $n$ sufficiently large.
\label{thm:asymptotic}
\end{theorem}

\subparagraph{Proof}
Applying Lemma \ref{lm:noncentral}, let $\lambda = \lambda_d = \lambda_k$.
Let
\[
\epsilon = 
P \left( \chi^2_{d-\ell}(\lambda) > c_{d,\alpha} \right) -
P \left( \chi^2_{k-\ell}(\lambda) > c_{k,\alpha} \right),
\]
which is strictly positive by virtue of Lemma {$\mathbf \chi^2$}.
Now use (\ref{eq:altd1}) and (\ref{eq:altk1}) to
choose $n \geq N$ sufficiently large that
\[
\left| P_{\vartheta + h/\sqrt{n}} 
\left( \Lambda_{d,n} > c_{d,\alpha} \right) -
P \left( \chi^2_{d-\ell} \left( \lambda_d \right) > c_{d,\alpha} \right)
\right| < \epsilon/2
\]
and
\[
\left| P_{\vartheta + h/\sqrt{n}} 
\left( \Lambda_{k,n} > c_{k,\alpha} \right) -
P \left( \chi^2_{k-\ell} \left( \lambda_k \right) > c_{k,\alpha} \right)
\right| < \epsilon/2.
\]
\hfill $\Box$

\bigskip
Under classical regularity conditions or, more generally,
under any set of conditions that ensure (\ref{eq:null})--(\ref{eq:altd2}),
Theorem \ref{thm:asymptotic} establishes that
the restricted LRT $\phi_d$ is asymptotically more powerful than the unrestricted LRT $\phi_k$ for local alternatives.
The following section compares restricted and unrestricted LRTs for finite samples.

\section{Hardy-Weinberg Equilibrium}
\label{HW}

Consider an experiment with three possible outcomes.
The model $\mbox{Trinomial}(\theta)$ specifies that the outcomes occur with probabilities $\theta=(\theta_1,\theta_2,\theta_3)$.
It is parametrized by the unit simplex in $\Re^3$,
$\Theta = \{ \theta \in [0,1]^3 : \theta_1+\theta_2+\theta_3=1 \}$.
  
Suppose that one draws $n$ i.i.d.\ observations from $\mbox{Trinomial}(\theta)$ and counts $x=(x_1,x_2,x_3)$, where $x_i$ records the number of occurrences of outcome $i$.  The unrestricted likelihood function of $\theta$ is
\[
L_x(\theta) =
P_\theta \left( X = \left( x_1,x_2,x_3 \right) \right) =
\frac{n!}{x_1! x_2! x_3!} \theta_1^{x_1} \theta_2^{x_2} \theta_3^{x_3},
\]
and the unrestricted maximum likelihood estimate of $\theta$ is
$\hat{\theta} = (x_1/n,x_2/n,x_3/n)$.
The unrestricted LRT rejects $H_0: \theta=\bar{\theta}$ if and only if
\[
\Lambda_2(x) = -2 \log \frac{L_x(\bar{\theta})}{L_x(\hat{\theta})} =
-2 \log \frac{\bar{\theta}^{x_1}\bar{\theta}^{x_2}\bar{\theta}^{x_3}}{
\hat{\theta}^{x_1}\hat{\theta}^{x_2}\hat{\theta}^{x_3}} 
\]
is sufficiently large.

Define $\psi : [0,1] \rightarrow \Theta$ by 
$\psi(\tau) = \left( \tau^2, 2\tau(1-\tau), (1-\tau)^2 \right)$.
The Hardy-Weinberg subfamily of trinomial distributions is 
parametrized by the embedded submanifold $\Psi = \{ \psi(\tau) : \tau \in [0,1] \}$.
Notice that $\mbox{dim } \Psi =1 < 2 = \mbox{dim } \Theta$.

Writing $\mbox{HW}(\tau) = \mbox{Trinomial}(\psi(\tau))$ and $m=2x_1+x_2$,
the likelihood function of $\tau$ is
\begin{eqnarray*}
L_x(\psi(\tau)) & = & 
P_{\psi(\tau)} \left( X = \left( x_1,x_2,x_3 \right) \right) =
\frac{n!}{x_1! x_2! x_3!} [\tau^2]^{x_1} [2\tau(1-\tau)]^{x_2} [(1-\tau)^2]^{x_3} \\ & = &
\frac{n!}{x_1! x_2! x_3!} 2^{x_2} \tau^m (1-\tau)^{2n-m},
\end{eqnarray*}
the maximum likelihood estimate of $\tau$ is $\hat{\tau} = m/(2n)$,
the restricted maximum likelihood estimate of $\theta$ is
$\tilde{\theta} = \psi(\hat{\tau})$, and
the restricted LRT rejects $H_0: \theta=\psi(\bar{\tau})$ if and only if
\[
\Lambda_1(x) = -2 \log \frac{L_x(\psi(\bar{\tau}))}{ L_x(\psi(\hat{\tau}))} =
-2 \log \frac{\bar{\tau}^m(1-\bar{\tau})^{2n-m}}{
\hat{\tau}^m(1-\hat{\tau})^{2n-m}} 
\]
is sufficiently large.

\subparagraph{Counterexample}
The trinomial experiment with $n=3$ has $10$ possible outcomes, enumerated in the first column of Table~\ref{tbl:counter}.
The three outcomes with the largest values of $\Lambda_2(x)$ are $(3,0,0)$, $(2,1,0)$, and $(2,0,1)$.  Denote this set of outcomes by $C_2$ and let 
\[
\alpha = P_{\psi(0.3)}(C_2) = 0.000729+0.010206+0.011907 = 0.022842,
\]
so that $C_2$ is the critical region for the unrestricted LRT of $H_0: \theta=\psi(0.3)$ of size $\alpha$.  In contrast,
the three outcomes with the largest values of $\Lambda_1(x)$ are $(3,0,0)$, $(2,1,0)$, and $(0,0,3)$.  Because 
\[
P_{\psi(0.3)} (0,0,3)  = 0.117649 > 0.011907 =
P_{\psi(0.3)} (2,0,1) ,
\]
the restricted LRT of $H_0: \theta=\psi(0.3)$ must be randomized in order to have size $\alpha$.  The randomized test will reject $H_0$ with certainty if $x=(3,0,0)$ or $x=(2,1,0)$, and with probability $0.011907/0.117649$ if $x=(0,0,3)$.

\begin{table}[tb]
\[
\begin{array}{|crrcrr|} \hline
x_1,x_2,x_3 & L_x(\psi(0.3)) & L_x(\hat{\theta}) & 
L_x(\psi(\hat{\tau})) & \Lambda_2(x) & \Lambda_1(x) \\ \hline
3,0,0 & 0.000729 & 27/27 &                            1 & 14.447674 & 14.447674 \\
2,1,0 & 0.010206 & 12/27 &  6 \cdot 5^5 \cdot 1^1 / 6^6 &  7.547699 &  7.346343 \\
2,0,1 & 0.011907 & 12/27 &  3 \cdot 4^4 \cdot 2^2 / 6^6 &  7.239397 &  3.420312 \\
1,2,0 & 0.047628 & 12/27 & 12 \cdot 4^4 \cdot 2^2 / 6^6 &  4.466808 &  3.420312 \\
1,1,1 & 0.111132 &  6/27 & 12 \cdot 3^3 \cdot 3^3 / 6^6 &  1.385918 &  1.046120 \\
1,0,2 & 0.064827 & 12/27 &  3 \cdot 2^2 \cdot 4^4 / 6^6 &  3.850206 &  0.031121 \\
0,3,0 & 0.074088 & 27/27 &  8 \cdot 3^3 \cdot 3^3 / 6^6 &  5.205003 &  1.046120 \\
0,2,1 & 0.259308 & 12/27 & 12 \cdot 2^2 \cdot 4^4 / 6^6 &  1.077617 &  0.031121 \\
0,1,2 & 0.302526 & 12/27 &  6 \cdot 1^1 \cdot 5^5 / 6^6 &  0.769316 &  0.567960 \\
0,0,3 & 0.117649 & 27/27 &                            1 &  4.280099 &  4.280099 \\ \hline
\end{array}
\]
\caption{Unrestricted (Trinomial) and restricted (Hardy-Weinberg) LRTs of $H_0:\theta=\psi(0.3)$ with $n=3$ observations.  Columns 1--2 list the possible outcomes and their exact probabilities under $H_0$; Column 3--4 list exact probabilities under the most likely Trinomial and Hardy-Weinberg distributions; Columns 5--6 list the unrestricted and restricted LRT statistics.}
\label{tbl:counter}
\end{table}

It is now apparent that the relative powers of the unrestricted and restricted LRTs at an alternative will depend on the probabilities of observing $(2,0,1)$ and $(0,0,3)$:  the restricted LRT will be more powerful at $\theta=\psi(\tau)$ if and only if
\begin{equation}
\frac{0.011907}{0.117649} P_{\psi(\tau)} (0,0,3) > P_{\psi(\tau)} (2,0,1).
\label{eq:cex}
\end{equation}
Some calculation reveals that (\ref{eq:cex}) holds when $\tau < 0.3$, but not when $\tau > 0.3$.  For $\tau \in (0.3,1)$,
the restricted LRT is {\em less}\/ powerful than the unrestricted LRT.

\hfill $\Box$

Our explication of the Counterexample reveals three key elements that allow construction of other counterexamples: first, a choice of $\bar{\tau}$ that causes $\Lambda_2$ and $\Lambda_1$ to order the possible outcomes differently; second, a level $\alpha$ for which the critical regions $C_2$ and $C_1$ differ with respect to a single outcome; and third, an alternative $\tau \neq \bar{\tau}$ under which the probability of $C_2$ exceeds the probability of $C_1$.
The remainder of this section is devoted to demonstrating that these elements occur far more generally.

In what follows, it will be convenient to work with the inverse likelihood ratios,
\[
R_2(x) = \frac{L_x \left( \hat{\theta} \right)}{
L_x \left( \psi \left( \bar{\tau} \right) \right)} =
\frac{ \left( \frac{x_1}{n} \right)^{x_1} \left( \frac{x_2}{n} \right)^{x_2} \left( \frac{x_3}{n} \right)^{x_3}}{2^{x_2} \bar{\tau}^m \left( 1-\bar{\tau} \right)^{2n-m}}
\]
and
\[
R_1(x) = \frac{L_x \left( \psi \left( \hat{\tau} \right) \right)}{
L_x \left( \psi \left( \bar{\tau} \right) \right)} =
\frac{ \left( \frac{m}{2n} \right)^m \left( 1-\frac{m}{2n} \right)^{2n-m}}{ \bar{\tau}^m \left( 1-\bar{\tau} \right)^{2n-m}},
\]
instead of $\Lambda_2$ and $\Lambda_1$.
Let $R_{2i}$ denote the level sets of $R_2$, ordered from largest $R_2$ value to smallest $R_2$ value, and let
\[
\alpha_{2j} = P_{\psi(\bar{\tau})} \left( \bigcup_{i=1}^j R_{2i} \right) = \sum_{i=1}^j P_{\psi(\bar{\tau})} \left( R_{2i} \right).
\]
For testing $H_0: \theta=\psi(\bar{\tau})$,
the unrestricted LRT of size $\alpha_{2j}$ has critical region
\[
C_2 \left( \alpha_{2j} \right) = \bigcup_{i=1}^j R_{2i},
\]
i.e., it rejects $H_0$ if and only if $x \in C_2(\alpha_{2j})$.
To obtain the LRT of size $\alpha \in (\alpha_{2j},\alpha_{2,j+1})$, it is necessary to randomize.  The conventional randomized LRT rejects $H_0$ with probability one if $x \in C_2(\alpha_{2j})$, with probability
$(\alpha-\alpha_{2j})/(\alpha_{2,j+1}-\alpha_{2j})$
if $x \in R_{2,j+1}$, and with probability zero otherwise.

The restricted case is analogous.  
Let $R_{1i}$ denote the level sets of $R_1$, ordered from largest $R_1$ value to smallest $R_1$ value, and let
\[
\alpha_{1j} = P_{\psi(\bar{\tau})} \left( \bigcup_{i=1}^j R_{1i} \right) = \sum_{i=1}^j P_{\psi(\bar{\tau})} \left( R_{1i} \right).
\]
For testing $H_0: \theta=\psi(\bar{\tau})$,
the restricted LRT of size $\alpha_{1j}$ has critical region
\[
C_1 \left( \alpha_{1j} \right) = \bigcup_{i=1}^j R_{1i},
\]
i.e., it rejects $H_0$ if and only if $x \in C_1(\alpha_{1j})$.
The conventional randomized LRT 
of size $\alpha \in (\alpha_{1j},\alpha_{1,j+1})$
rejects $H_0$ with probability one if $x \in C_1(\alpha_{1j})$, with probability
$(\alpha-\alpha_{1j})/(\alpha_{1,j+1}-\alpha_{1j})$
if $x \in R_{1,j+1}$, and with probability zero otherwise.

Notice that $R_1$ is constant on the level sets of the integer-valued random variable $M=2X_1+X_2$.
For $R_1$ to assume the same value with $m_1 \neq m_2$, it must be that
\[
\frac{ \left( \frac{m_1}{2n} \right)^{m_1} \left( 1-\frac{m_1}{2n} \right)^{2n-m_1}}{ \bar{\tau}^{m_1} \left( 1-\bar{\tau} \right)^{2n-m_1}} =
\frac{ \left( \frac{m_2}{2n} \right)^{m_2} \left( 1-\frac{m_2}{2n} \right)^{2n-m_2}}{ \bar{\tau}^{m_2} \left( 1-\bar{\tau} \right)^{2n-m_2}},
\]
which is equivalent to
\begin{equation}
\left( \frac{\bar{\tau}}{1-\bar{\tau}} \right)^{m_1-m_2} =
\frac{m_1^{m_1} \left( 2n-m_1 \right)^{2n-m_1}}{
m_2^{m_2} \left( 2n-m_2 \right)^{2n-m_2}}.
\label{eq:rational}
\end{equation}
Because the right-hand side of (\ref{eq:rational}) is rational,
(\ref{eq:rational}) cannot obtain if $\bar{\tau}$ is irrational.
We thus establish

\begin{lemma}
If $\bar{\tau}$ is irrational, then the level sets of $R_1$ coincide with the level sets of $2X_1+X_2$.  In particular, each level set of $R_1$ is associated with a single value of $2X_1+X_2$.
\label{lm:irrational1}
\end{lemma}

The level sets of $R_2$ are not so easily characterized, but suppose
that $R_2(x)=R_2(y)$.  Set $m_1=2x_1+x_2$ and $m_2=2y_1+y_2$.
Then
\[
\frac{\bar{\tau}^{m_1} (1-\bar{\tau})^{2n-m_1}}{
\bar{\tau}^{m_2} (1-\bar{\tau})^{2n-m_2}} =
\frac{(2x_1)^{x_1}x_2^{x_2}(2x_3)^{x_3}}{
(2y_1)^{y_1}y_2^{y_2}(2y_3)^{y_3}},
\]
the right-hand side of which is rational and
the left-hand side of which is irrational if $\bar{\tau}$ is irrational and $m_1 \neq m_2$.  We thus establish

\begin{lemma}
If $\bar{\tau}$ is irrational, then each level set of $R_2$ is associated with a single value of $2X_1+X_2$.
\label{lm:irrational2}
\end{lemma}

Next we compare how $R_1$ and $R_2$ arrange possible outcomes into critical regions.  Our first result states that the restricted and unrestricted LRTs agree on which outcome (or set of outcomes) is most adverse to $H_0:\theta=\psi(\bar{\tau})$.

\begin{lemma}
For every $\bar{\tau} \in (0,1)$, $R_{11}=R_{21}$.
\label{lm:worst}
\end{lemma}
\subparagraph{Proof}
The level set $R_{11}$ consists of the outcomes that maximize
\[
R_1(x) = \frac{ \left( \frac{m}{2n} \right)^m \left( 1-\frac{m}{2n} \right)^{2n-m}}{ \bar{\tau}^m \left( 1-\bar{\tau} \right)^{2n-m}}.
\]
If $\bar{\tau} \in (0,0.5)$,
then the denominator is minimized by $m=2n$,
which also maximizes the numerator.
In this case, $R_{11} = \{ (n,0,0) \}$.
By the same reasoning, if $\bar{\tau} \in (0.5,1)$,
then $R_{11} = \{ (0,0,n) \}$.
If $\bar{\tau}=0.5$, then the denominator is constant and the numerator is maximized by either $m=2n$ or $m=0$; hence,
$R_{11} = \{ (n,0,0),(0,0,n) \}$.

The level set $R_{21}$ consists of the outcomes that maximize
\[
R_2(x) = 
\frac{ \left( \frac{x_1}{n} \right)^{x_1} \left( \frac{x_2}{n} \right)^{x_2} \left( \frac{x_3}{n} \right)^{x_3}}{2^{x_2} \bar{\tau}^m \left( 1-\bar{\tau} \right)^{2n-m}}.
\]
If $\bar{\tau} \in (0,0.5)$,
then the denominator is minimized by $(n,0,0)$,
which also maximizes the numerator.
Hence, $R_{21} = \{ (n,0,0) \} = R_{11}$.
By the same reasoning, if $\bar{\tau} \in (0.5,1)$,
then $R_{21} = \{ (n,0,0) \} = R_{11}$.
If $\bar{\tau}=0.5$, then the denominator is minimized by any $x$ with $x_2=0$ and the numerator is maximized by either
$(n,0,0)$ or $(0,0,n)$; hence,
$R_{21} = \{ (n,0,0),(0,0,n) \} = R_{11}$.

\hfill $\Box$

Next we establish some conditions under which $R_2$ and $R_1$ induce different orderings of the possible outcomes.  More precisely, the conditions in Lemma~\ref{lm:order} ensure the
existence of outcomes $x$ and $y$ such that $R_2(x) > R_2(y)$ and $R_1(x) < R_1(y)$.  Note the strict inequalities in this definition, which are essential to the proof of Theorem~\ref{thm:HWirrational}.

\begin{lemma}
The unrestricted and restricted LRTs order the possible outcomes differently under any of the following conditions:
\begin{enumerate}
\item $\bar{\tau} \in (0.2,1/3) \cup (2/3,0.8)$;
\item $\bar{\tau} \in (1/3,2/3)$ and $n \geq 2$;
\item $\bar{\tau} \in \{ 1/3, 2/3 \}$ and $n \geq 3$.
\end{enumerate}
\label{lm:order}
\end{lemma}
\subparagraph{Proof}
First, consider the three extremal outcomes displayed in Table~\ref{tbl:extremal0}.  The ordering of these outcomes is determined by $R_2(x)=L_x(\hat{\theta})/L_x(\psi(\bar{\tau}))$ for the unrestricted LRT and by $R_1(x)=L_x(\psi(\hat{\tau}))/L_x(\psi(\bar{\tau}))$ for the restricted LRT.  For $\bar{\tau} \in (0.2,1/3)$,
\[
\frac{1}{2 \bar{\tau} (1-\bar{\tau})} >
\frac{1}{(1-\bar{\tau})^2}
\]
and the unrestricted LRT orders $(0,n,0) \succ (0,0,n)$, whereas
\[
\frac{1}{(1-\bar{\tau})^2} >
\frac{1/2}{2 \bar{\tau} (1-\bar{\tau})}
\]
and the restricted LRT orders $(0,0,n) \succ (0,n,0)$.
By symmetry, the unrestricted and restricted orders also differ if $\bar{\tau} \in (2/3,0.8)$.

\begin{table}[tb]
\[
\begin{array}{|cccc|} \hline
x_1,x_2,x_3 & L_x(\psi(\bar{\tau})) & L_x(\hat{\theta}) & 
L_x(\psi(\hat{\tau})) \\ \hline
(n,0,0) & \bar{\tau}^{2n} & 1 & 1 \\
(0,n,0) & 2^n \bar{\tau}^n (1-\bar{\tau})^n & 1 & 2^{-n} \\
(0,0,n) & (1-\bar{\tau})^{2n} & 1 & 1 \\ \hline
\end{array}
\]
\caption{Probabilities of three extremal outcomes under
$\theta=\psi(\bar{\tau})$ (null hypothesis),
$\theta=\hat{\theta}$ (unrestricted MLE), and
$\theta=\psi(\hat{\tau})$ (restricted MLE).}
\label{tbl:extremal0}
\end{table}

Next, assume that $n \geq 2$ and consider the outcomes displayed in Table~\ref{tbl:extremal1}.  The unrestricted LRT places $(n-1,1,0) \succ (n-1,0,1)$ if and only if 
\[
2 \bar{\tau}^{2n-1} (1-\bar{\tau}) <
\bar{\tau}^{2n-2} (1-\bar{\tau})^{2},
\]
which obtains if and only if $\bar{\tau} < 1/3$.
In contrast, the restricted LRT places $(n-1,1,0) \succ (n-1,0,1)$ if and only if 
\[
\left( \frac{2n-1}{2n} \right)^{2n-1} \frac{1}{2n} \div
\bar{\tau}^{2n-1} (1-\bar{\tau}) >
\left( \frac{2n-2}{2n} \right)^{2n-2} 
\left( \frac{2}{2n} \right)^{2} \div
\bar{\tau}^{2n-2} (1-\bar{\tau})^{2},
\]
which obtains if and only if
\[
\bar{\tau} < b(n) =
\frac{(2n-1)^{2n-1}}{(2n-1)^{2n-1}+4(2n-2)^{2n-2}}.
\]
If $n=2$, then
\[
4(2n-2)^{2n-2} = 16<27 = (2n-1)^{2n-1}
\]
and $b(n)>0.5$.
If $n \geq 3$, then
\[
4(2n-2)^{2n-2} \leq (2n-2)^{2n-1} < (2n-1)^{2n-1}
\]
and again $b(n) > 0.5$.  It follows
that the unrestricted and restricted LRT orders differ if $\bar{\tau} \in (1/3,0.5]$.  By symmetry, they also differ if
$\bar{\tau} \in [0.5,2/3)$.

\begin{table}[tb]
\[
\begin{array}{|cccc|} \hline
x_1,x_2,x_3 & L_x(\psi(\bar{\tau})) & L_x(\hat{\theta}) & 
L_x(\psi(\hat{\tau})) \\ \hline
(n-1,1,0) & 2n \bar{\tau}^{2n-1} (1-\bar{\tau}) &
\left( \frac{n-1}{n} \right)^{n-1} &
2n \left( \frac{2n-1}{2n} \right)^{2n-1} 
\frac{1}{2n} \\
(n-1,0,1) & n \bar{\tau}^{2n-2} (1-\bar{\tau})^{2} &
\left( \frac{n-1}{n} \right)^{n-1} &
n \left( \frac{2n-2}{2n} \right)^{2n-2} 
\left( \frac{2}{2n} \right)^{2} \\
(1,n-1,0) & 2^{n-1}n \bar{\tau}^{n+1} (1-\bar{\tau})^{n-1} &
\left( \frac{n-1}{n} \right)^{n-1} &
2^{n-1}n \left( \frac{n+1}{2n} \right)^{n+1} 
\left( \frac{n-1}{2n} \right)^{n-1} \\
(0,n-1,1) & 2^{n-1}n \bar{\tau}^{n-1} (1-\bar{\tau})^{n+1} &
\left( \frac{n-1}{n} \right)^{n-1} &
2^{n-1}n \left( \frac{n-1}{2n} \right)^{n-1} 
\left( \frac{n+1}{2n} \right)^{n+1} \\
(1,0,n-1) & n \bar{\tau}^{2} (1-\bar{\tau})^{2n-2} &
\left( \frac{n-1}{n} \right)^{n-1} &
n \left( \frac{2}{2n} \right)^{2} 
\left( \frac{2n-2}{2n} \right)^{2n-2} \\
(0,1,n-1) & 2n \bar{\tau} (1-\bar{\tau})^{2n-1} &
\left( \frac{n-1}{n} \right)^{n-1} &
2n \frac{1}{2n} 
\left( \frac{2n-1}{2n} \right)^{2n-1} \\ \hline
\end{array}
\]
\caption{Probabilities of six more outcomes under
$\theta=\psi(\bar{\tau})$ (null hypothesis),
$\theta=\hat{\theta}$ (unrestricted MLE), and
$\theta=\psi(\hat{\tau})$ (restricted MLE), assuming $n>1$.}
\label{tbl:extremal1}
\end{table}

Finally, let $\bar{\tau}=1/3$ and $n \geq 3$.  Then
\[
n \left( \frac{1}{3} \right)^2 \left( \frac{2}{3} \right)^{2n-2} =
n \frac{2^{n-2}}{3^{2n}} < n \frac{2^{2n}}{3^{2n}} =
2n \frac{1}{3} \left( \frac{2}{3} \right)^{2n-1}
\]
and it follows that the unrestricted LRT orders
$(1,0,n-1) \succ (0,1,n-1)$.
In contrast,
\[
\frac{L_{(1,0,n-1)}((\psi(\hat{\tau}))}{L_{(1,0,n-1)}(\psi(1/3))} \div
\frac{L_{(0,1,n-1)}((\psi(\hat{\tau}))}{L_{(0,1,n-1)}(\psi(1/3))} =
c(n) = \frac{8(2n-2)^{2n-2}}{(2n-1)^{2n-1}}.
\]
If $n=3$, then
\[
8(2n-2)^{2n-2} = 2048 < 3125 = (2n-1)^{2n-1}
\]
and $c(n)<1$.
If $n=4$, then
\[
8(2n-2)^{2n-2} = 373248 < 823543 = (2n-1)^{2n-1}
\]
and $c(n)<1$.
If $n \geq 5$, then
\[
8(2n-2)^{2n-2} \leq (2n-2)^{2n-1} < (2n-1)^{2n-1}
\]
and $c(n)<1$.
It follows that the restricted LRT orders
$(1,0,n-1) \prec (0,1,n-1)$.

By symmetry, the unrestricted and restricted LRT orders differ if $\bar{\tau}=2/3$ and $n \geq 3$.

\hfill $\Box$

We now establish our crucial result.

\begin{theorem}
Let $\psi$ parametrize the Hardy-Weinberg submodel of the trinomial experiment.
Suppose that $\bar{\tau} \in (0,1)$ is irrational, 
and that $(\bar{\tau},n)$ is such that the restricted and unrestricted LRTs of $H_0:\theta=\psi(\bar{\tau})$ order the possible outcomes of the trinomial experiment differently.
Then there exist $\alpha \in (0,1)$ and $\tau \in (0,1)$ such that the restricted LRT of size $\alpha$ is less powerful at $\tau$ than the unrestricted LRT of size $\alpha$.
\label{thm:HWirrational}
\end{theorem} 
\subparagraph{Proof}
Lemma~\ref{lm:worst} states that $R_{11}=R_{21}$.
If $R_{1i}=R_{2i}$ for $i=1,\ldots,j$ and
\[
x,y \in \bigcup_{i=1}^j R_{1i} = \bigcup_{i=1}^j R_{2i},
\]
then it cannot be that 
$R_2(x) > R_2(y)$ and $R_1(x) < R_1(y)$.
Hence, for the restricted and unrestricted LRTs to order the possible outcomes differently, there must exist a value of $j$ for which $R_{1j} \neq R_{2j}$
Let $j*$ denote the smallest such value of $j$.
Because of Lemmas \ref{lm:irrational1} and \ref{lm:irrational2}, there are two possibilities: either
\begin{enumerate}

\item $R_{1j*}$ and $R_{2j*}$ are associated with different values of $2X_1+X_2$; or

\item $R_{1j*}$ and $R_{2j*}$ are associated with the same value of $2X_1+X_2$, in which case $R_{1j*} \neq R_{2j*}$ implies that $R_{2j*}$ is a proper subset of $R_{1j*}$.

\end{enumerate}

The first case is straightforward.
Let $m_1$ and $m_2$ denote the values of $2X_1+X_2$ associated with 
$R_{1j*}$ and $R_{2j*}$.
Let $\alpha = \min(\alpha_{1j*},\alpha_{2j*})$,
so that the restricted and unrestricted LRTs of size $\alpha$ have critical regions that are identical except for (possibly randomized) outcomes in $R_{1j*}$ or $R_{2j*}$.  Any power differences between the two LRTs will accrue from these outcomes.

The probability that the restricted LRT will reject $H_0$ as a result of 
$x \in R_{1j*}$ is
\[
r_1 P_{\psi(\tau)} \left( R_{1i*} \right) =
r_1 n! \tau^{m_1} (1-\tau)^{2n-m_1} \sum_{x \in R_{1i*}} 
\frac{2^{x_2}}{x_1!x_2!x_3!},
\]
and the probability that the unrestricted LRT will reject $H_0$ as a result of 
$x \in R_{2j*}$ is
\[
r_2 P_{\psi(\tau)} \left( R_{2j*} \right) =
r_2 n! \tau^{m_2} (1-\tau)^{2n-m_2} \sum_{x \in R_{2j*}}
\frac{2^{x_2}}{x_1!x_2!x_3!},
\]
where $r_1$ and $r_2$ are randomization factors.
If $\tau=\bar{\tau}$, then these quantities are equal.
If
\[
r_1 P_{\psi(\tau)} \left( R_{1j*} \right) <
r_2 P_{\psi(\tau)} \left( R_{2j*} \right),
\]
which obtains when
\[
\left( \frac{\tau}{1-\tau} \right)^{m_1-m_2} <
 \left( r_2 \sum_{x \in R_{2j*}} 
 \frac{2^{x_2}}{x_1!x_2!x_3!} \right) \div
 \left( r_1 \sum_{x \in R_{1j*}} 
 \frac{2^{x_2}}{x_1!x_2!x_3!} \right),
\]
then the unrestricted LRT has greater power at $\tau$ then the restricted LRT.

The key to the preceding argument lies in identifying a size for which the critical regions of the respective LRTs differ only in replacing one level set of $2X_1+X_2$ with another of different value.
If $R_{1j*}$ and $R_{2j*}$ are associated with the same value of $2X_1+X_2$, then identifying such a size is more complicated.
If $R_{2j*}$ is a proper subset of $R_{1j*}$, 
then we set $\alpha = \min(\alpha_{1j*},\alpha_{2,j*+1})$
and consider $R_{2,j*+1}$.
If $R_{2,j*+1}$ is associated with a different value of $2X_1+X_2$,
then the probability that the unrestricted LRT will reject $H_0$ as a result of 
$x \in R_{2j*} \cup R_{2,j*+1}$ is
\begin{eqnarray*}
P_{\psi(\tau)} \left( R_{2j*} \right) + r_2
P_{\psi(\tau)} \left( R_{2,j*+1} \right) & = &
n! \tau^{m_1} (1-\tau)^{2n-m_1} \sum_{x \in R_{2j*}}
\frac{2^{x_2}}{x_1!x_2!x_3!} + \\
 & & r_2 n!
\tau^{m_2} (1-\tau)^{2n-m_2} \sum_{x \in R_{2,j*+1}}
\frac{2^{x_2}}{x_1!x_2!x_3!}.
\end{eqnarray*}
If
\[
r_1 P_{\psi(\tau)} \left( R_{1j*} \right) <
P_{\psi(\tau)} \left( R_{2j*} \right) + r_2
P_{\psi(\tau)} \left( R_{2,j*+1} \right),
\]
which obtains when
\begin{eqnarray*}
\left( \frac{\tau}{1-\tau} \right)^{m_1-m_2} & < &
 \left( r_2 \sum_{x \in R_{2,j*+1}}
 \frac{2^{x_2}}{x_1!x_2!x_3!} \right) \div \\ & &
 \left( r_1 \sum_{x \in R_{1j*}}
 \frac{2^{x_2}}{x_1!x_2!x_3!} -
\sum_{x \in R_{2j*}}
\frac{2^{x_2}}{x_1!x_2!x_3!} \right),
\end{eqnarray*}
then the unrestricted LRT has greater power at $\tau$ than the restricted LRT.

Continuing in this manner,
if $R_{2,j*+1}$ is not associated with a different value of $2X_1+X_2$, then we increment $i$ in 
$\alpha = \min(\alpha_{1j*},\alpha_{2,j*+i})$
and $R_{2,j*+i}$ until either we do encounter a different value or until both LRTs have the same critical region for size $\alpha_{1j*}=\alpha_{2,j*+i*}$.  
If we encounter a different value, then the same reasoning used in the previous paragraph establishes values of $\tau$ for which the unrestricted LRT has greater power at $\tau$ then the restricted LRT.

If we exhaust the outcomes in $R_{1j*}$ without encountering a different value, then we obtain
\begin{eqnarray*}
C_1 \left( \alpha_{1j*} \right) & = &
R_{11} \cup \cdots \cup R_{1,j*-1} \cup R_{1j*} \\ 
 & = &
R_{21} \cup \cdots \cup R_{2,j*-1} \cup \bigcup_{i=1}^{i*} R_{2,j*+i} \\
 & = & C_2 \left( \alpha_{1j*} \right).
\end{eqnarray*}
In this case, however, we have not yet discovered   
\[
x,y \in C_1 \left( \alpha_{1j*} \right) =
C_2 \left( \alpha_{1j*} \right)
\]
for which
$R_2(x) > R_2(y)$ and $R_1(x) < R_1(y)$.
We have already observed that no such reversal is possible for
\[
x,y \in \bigcup_{i=1}^{j*-1} R_{1i} = \bigcup_{i=1}^{j*-1} R_{2i},
\]
and it is certainly not possible for
\[
x,y \in R_{1j*} = \bigcup_{i=1}^{i*} R_{2,j*+i}
\]
because $R_{1j*}$ is a single level set of $R_1$.

If we reach the case of
$C_1 (\alpha_{1j*}) = C_2 (\alpha_{1j*})$, 
then we progress to $R_{1,j*+1}$ versus $R_{2,j*+i*+1}$ and apply the same reasoning.  Because there exist $x$ and $y$ for which
$R_2(x) > R_2(y)$ and $R_1(x) < R_1(y)$, we will eventually identify a size for which the critical regions of the two LRTs differ only in replacing one level set of $2X_1+X_2$ with another of different value, and we have already demonstrated how to derive the desired result when that circumstance obtains.

\hfill $\Box$

Combining Theorem~\ref{thm:HWirrational} and Lemma~\ref{lm:order}, we obtain the following result.

\begin{corollary}
Let $\psi$ parametrize the Hardy-Weinberg submodel of the trinomial experiment with $n \geq 2$ and consider the null hypothesis $H_0: \theta = \psi(\bar{\tau})$.
For almost every $\bar{\tau} \in (0.2,0.8)$, 
there exist $\alpha \in (0,1)$ and $\tau \in (0,1)$ such that the restricted LRT of size $\alpha$ is less powerful at $\tau$ than the unrestricted LRT of size $\alpha$.
\label{cor:HWae}
\end{corollary}

\section{Discussion}
\label{disc}

Marden's 1982 conjecture,
that restricting the set of alternatives increases the power of a likelihood ratio test, was falsified in 2003.  Unfortunately,
the counterexample constructed in \citep{abudayyeh&etal:2003}
is highly technical and involves a restriction that maintains the dimension of the alternative.
It is inviting to modify the general LRT conjecture and speculate that it holds if the submodel has lower dimension than the model.
As demonstrated in Section~\ref{asymptotic}, 
the dimension-restricted LRT conjecture does hold asymptotically.
However, as demonstrated in Section~\ref{HW},
the Hardy-Weinberg submodel of the trinomial experiment provides a wealth of counterexamples when sample size is finite.
This finding may surprise many readers, as it surprised us.

\begin{figure}[tb]
\begin{center}
\includegraphics[width=0.8\textwidth]{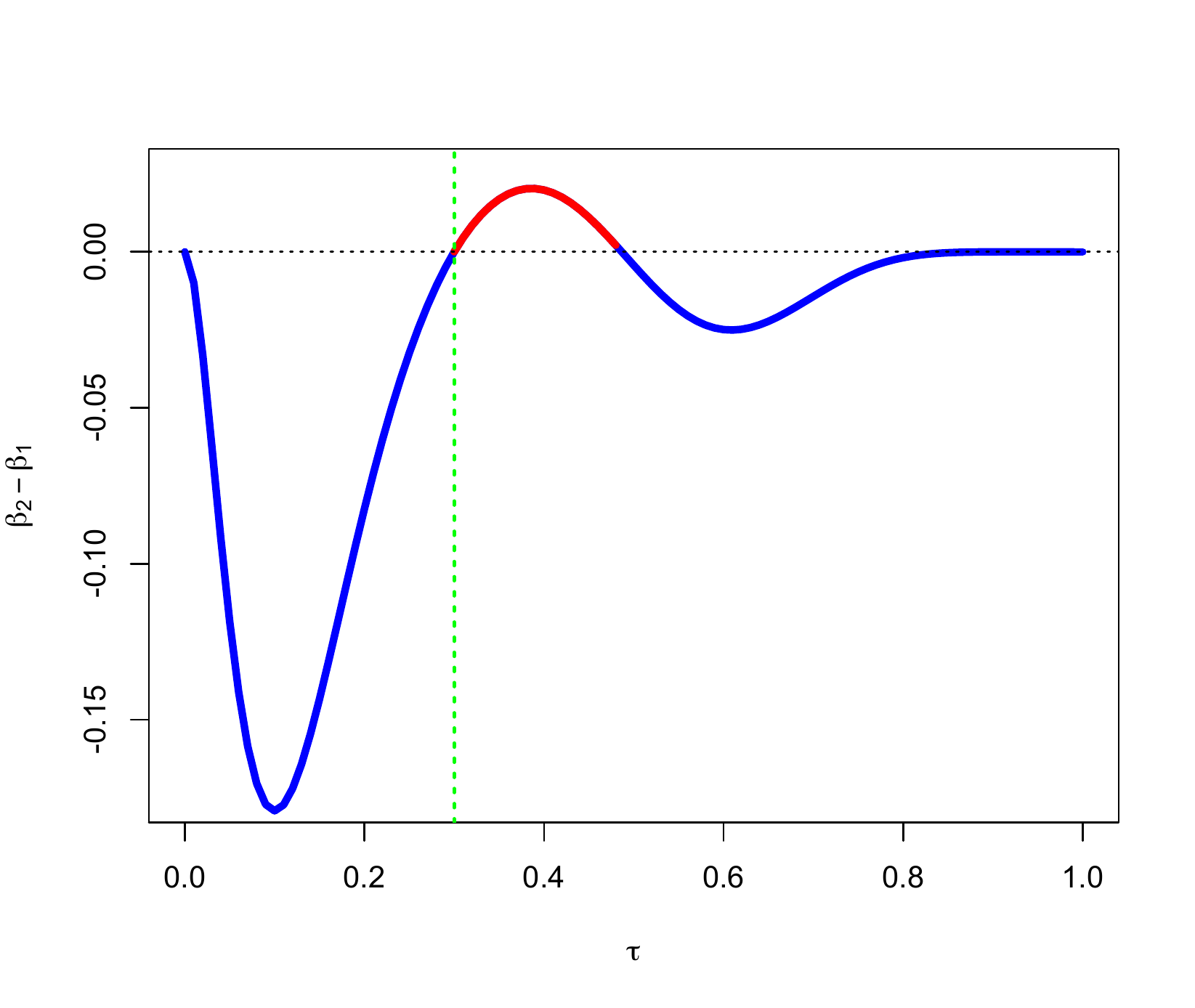} 
\end{center}
\caption{Power of the unrestricted LRT ($\beta_2$) minus power of the restricted LRT ($\beta_1$) for testing $H_0: \theta = \psi(0.3)$ with $\alpha=0.05$ and $n=10$.  The alternatives $\{ \theta=\psi(\tau) : \tau \in (0,1) \}$ are displayed on the horizontal axis.  The restricted LRT is superior where the power difference is less than zero (blue), inferior where it exceeds zero (red).} 
\label{fig:HWpower}
\end{figure}

Several observations are in order.  First, the counterexamples presented in Section~\ref{HW} only suggest that the restricted LRT is not {\em uniformly}\/ more powerful than the unrestricted LRT.  The fact that there exist sizes and alternatives for which the restricted LRT is less powerful does not imply that we should prefer the unrestricted LRT.
Our first counterexample, in which the restricted LRT is more powerful than the unrestricted LRT for alternatives in $(0,0.3)$ but less powerful for alternatives in $(0.3,1)$, is dramatic.
Nevertheless, in most of the Hardy-Weinberg examples that we have studied,
the restricted LRT outperforms the unrestricted LRT in more cases---and by wider margins---than the reverse.
We display an example in Figure~\ref{fig:HWpower}.

The sufficient conditions identified in Section~\ref{HW},
e.g., in Corollary~\ref{cor:HWae}, are expedient but hardly necessary.  We expect that more elaborate calculations will extend the scope of Lemma~\ref{lm:order}, while continuity arguments will extend the results of Theorem~\ref{thm:HWirrational} from irrational $\bar{\tau}$ to real intervals of $\bar{\tau}$.
This is not our present concern.

Because the Hardy-Weinberg subfamily of trinomial distributions is widely used to model genetic equilibrium, it provides a compelling counterexample to the dimension-restricted LRT conjecture. 
It is hardly unique, however, for we have constructed other submodels of multinomial models (including submodels of dimension and/or co-dimension greater than one) that also falsify the conjecture.
Our current efforts are focussed on identifying conditions under which the dimension-restricted LRT conjecture does---or does not---hold.

Our counterexample is also compelling from the perspective of statistical theory.  As a submodel of the trinomial model, the Hardy-Weinberg distributions constitute a curved exponential family.  Furthermore, as demonstrated in
\citep[Example~2.3.4]{kass&voss:1997},
these distributions can be parametrized as a regular $1$-parameter exponential family.
For testing simple null hypotheses, 
it follows from results in
\citep[Section~4.2]{lehmann:tsh2}
that there exists a uniformly most powerful test among all unbiased tests.  The difficulty is not that one cannot improve power by restricting the trinomial model to the Hardy-Weinberg submodel, but that LRTs may fail to do so.

For testing simple null and alternative hypotheses,
the Neyman-Pearson Lemma states that a test based on the ratio of the null and alternative probability densities is most powerful.
Likelihood ratio tests are heuristic extensions of this construction that possess various desirable asymptotic properties.
However, although LRTs are asymptotically unbiased
\citep[Section~53.3]{borokov:1999},
they may be biased for finite sample sizes.
The Counterexample in Section~\ref{HW}
considered the restricted LRT of $H_0: \psi(0.3)$ with size $\alpha = 0.22842$.  The power of that test at alternative $\tau$ is
\[
\tau^6 + 3 \cdot 2 \cdot \tau^5(1-\tau) +
\frac{0.011907}{0.117649} (1-\tau)^6,
\]
which numerical calculation reveals is slightly less than $\alpha$ for a small interval of alternatives. 

From a modeling perspective, our report is a cautionary tale: restricting inferences about a submodel to that submodel will strike most statisticians as natural and intuitive, yet doing so is not necessarily most powerful. As noted above, this finding is not due to some pathology unique to the Hardy-Weinberg submodel, which is nicely behaved; we have constructed additional counterexamples using various multinomial submodels, including submodels of dimension and/or co-dimension greater than one. Taken together, these examples constitute a compelling demonstration that efforts to exploit submodel structure may require considerable care.

\section*{Acknowledgments}
This work was partially supported by a National Security Science and Engineering Faculty Fellowship,
DARPA XDATA contract FA8750-12-2-0303,
DARPA SIMPLEX contract N66001-15-C-4041,
DARPA GRAPHS contract N66001-14-1-4028,
and NSF award DBI-1451081.

\bibliography{$HOME/lib/tex/stat,$HOME/lib/tex/mds,$HOME/lib/tex/math}

\end{document}